\documentclass[]{interact}

\usepackage{epstopdf}
\usepackage[caption=false]{subfig}

\usepackage{indentfirst}
\usepackage{url}
\usepackage{color}
\usepackage{float}
\usepackage{graphicx}   
\usepackage{hyperref}
\hypersetup{
    colorlinks=true,
    linkcolor=black,
    filecolor=magenta,      
    urlcolor=blue,
    citecolor=cyan,
}

\usepackage[numbers,sort&compress]{natbib}
\usepackage{amsmath}

\DeclareMathOperator*{\argmin}{arg\,min}
\newcommand\norm[1]{\left\lVert#1\right\rVert}
\bibpunct[, ]{[}{]}{,}{n}{,}{,}
\renewcommand\bibfont{\fontsize{10}{12}\selectfont}
\makeatletter
\def\NAT@def@citea{\def\@citea{\NAT@separator}}
\makeatother

\theoremstyle{plain}
\newtheorem{theorem}{Theorem}[section]

\theoremstyle{definition}

\theoremstyle{remark}
\newtheorem{remark}{Remark}

\begin{document}

\articletype{ARTICLE TEMPLATE}

\title{Angle-Aware Coverage with Camera Rotational Motion Control}

\author{
\name{Zhiyuan Lu\textsuperscript{a}\thanks{zhiyuan@hfg.sc.e.titech.ac.jp},
Muhammad Hanif\textsuperscript{\href{https://orcid.org/0009-0001-0475-8711}{\includegraphics[scale=0.2]{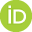}}a},
Takumi Shimizu\textsuperscript{\href{https://orcid.org/0000-0001-7604-5414}{\includegraphics[scale=0.2]{orcid_32x32.png}}a},
and Takeshi Hatanaka\textsuperscript{\href{https://orcid.org/0000-0003-3143-121X}{\includegraphics[scale=0.2]{orcid_32x32.png}}a}}
\affil{\textsuperscript{a}Tokyo Institute of Technology, 2-12-1 Ookayama, Meguro-ku, Tokyo, 152-8550, Japan}
}

\maketitle

\begin{abstract}

This paper presents a novel control strategy for drone networks to improve the quality of 3D structures reconstructed from aerial images by drones.  Unlike the existing coverage control strategies for this purpose, our proposed approach simultaneously controls both the camera orientation and drone translational motion, enabling more comprehensive perspectives and enhancing the map's overall quality. Subsequently, we present a novel problem formulation, including a new performance function to evaluate the drone positions and camera orientations. We then design a QP-based controller with a control barrier-like function for a constraint on the decay rate of the objective function. The present problem formulation poses a new challenge, requiring significantly greater computational efforts than the case involving only translational motion control. We approach this issue technologically, namely by introducing JAX, utilizing just-in-time (JIT) compilation and Graphical Processing Unit (GPU) acceleration. We finally conduct extensive verifications through simulation in ROS (Robot Operating System) and show the real-time feasibility of the controller and the superiority of the present controller to the conventional method.

\end{abstract}


\section{Introduction}


The advancement of 3D map reconstruction technology has played a pivotal role in supporting various sectors, including building information modeling \cite{zheng2018multi,macher2017point}, precision agriculture \cite{comba2018unsupervised,dong20174d}, and construction site inspection \cite{xue2021review}. The success of the Structure from Motion (SfM) algorithm \cite{schonberger2016structure} has been instrumental in efficiently reconstructing a 3D model of a target object by analyzing a collection of images captured by a camera. Notably, the past decades have witnessed the introduction of numerous sensors and platforms for capturing high-quality images, with unmanned aerial vehicles (UAVs) emerging as particularly prevalent platforms capable of autonomous image acquisition \cite{maboudi2023review}. Therefore, ensuring the efficient coverage of viewpoints by UAVs becomes a crucial requirement for achieving high-quality map reconstruction using the SfM algorithm.


To ensure efficient coverage of viewpoints, at least two crucial factors must be considered. Firstly, achieving comprehensive observation of every part of the scene from multiple angles is vital to enhance the reliability and accuracy of 3D reconstruction \cite{maboudi2023review}. This necessitates sufficient overlap between viewpoints, compelling UAVs to explore a 6D configuration encompassing both 3D position and 3D angle orientation. However, in practice, the roll angle is often neglected, simplifying the search space to a 5D configuration for camera poses \cite{liu2021aerial,smith2018aerial}. Secondly, addressing the one-time visit problem is essential \cite{shimizu2021angle}. This involves systematically capturing images across the target field, ensuring that each 5D point within the target field is observed and visited.


Addressing the aforementioned requirement, classical coverage control algorithms have traditionally been widely adopted to tackle the scenario outlined above \cite{cortes2005spatially}. Moreover, various studies have also explored coverage control strategies in environments with diverse shapes \cite{kantaros2014visibility,kantaros2015distributed,schwager2006distributed,schwager2009decentralized,breitenmoser2010voronoi,stergiopoulos2015distributed}. However, a notable limitation in the existing literature is that they do not consider coverage in a 5D search space, a crucial aspect in the context of UAV-based 3D map reconstruction. Furthermore, many of these studies tend to guide robots into static configurations, falling short of covering each point within the target field.


Alternatively, there has been a study on persistent coverage control algorithms, enabling UAVs to monitor the environment persistently \cite{palacios2016distributed, sugimoto2015experimental, dan2020control, dan2021persistent}. However, once more, these studies predominantly focus on 2D coverage, overlooking the necessity to observe points from multiple angles. Furthermore, persistent coverage strategies may not address the specific requirement of a one-time visit, as specified earlier.


Recent work in angle-aware coverage control \cite{shimizu2021angle} marks the initial attempt at achieving 5D coverage with multiple UAVs. This approach also provides a solution to the one-time visit problem. By integrating the concept of capturing images from various angles within the coverage control framework, significant improvements in 3D map reconstruction quality have been demonstrated, as shown in \cite{suenaga2022experimental}. However, it is essential to note that even in \cite{shimizu2021angle}, the drones' camera orientations remained fixed, suggesting room for further improvements. This could be achieved by considering dynamic camera control of the UAV using gimbal mechanisms. In an earlier exploration into camera angle control, documented in a previous study \cite{schwager2011eyes,hatanaka2019visual}, the focus was primarily on adjusting the camera's angle while keeping the camera position static; therefore, it could not solve the one-time visit problem.


In the context of this paper, our objective is to extend the existing angle-aware coverage control framework \cite{shimizu2021angle} by integrating camera rotational motion control via gimbal mechanisms. This enhancement allows for simultaneous control of both camera orientation and drone motion within the coverage control framework. To achieve this, we present a new formulation of the problem that includes a new performance function to assess the camera orientations and drone positions. We then design a QP-based controller \cite{ames2016control} with a constraint using a control barrier-like function on the decay rate of the objective function. This new problem formulation is significantly more computationally demanding than the case of coverage with translational motion control only. We address this challenge by implementing JAX \cite{jax2018github} , employing JIT compilation and GPU acceleration. Finally, using simulation in the ROS, we conduct thorough verifications demonstrating the controller's real-time viability and superiority over the traditional approach.

\begin{figure}[t]\
    \centering
    \subfloat[3D Map Reconstruction Scene]{
        \label{fig:3d_map_recon_scene}
        \includegraphics[width = 13cm]{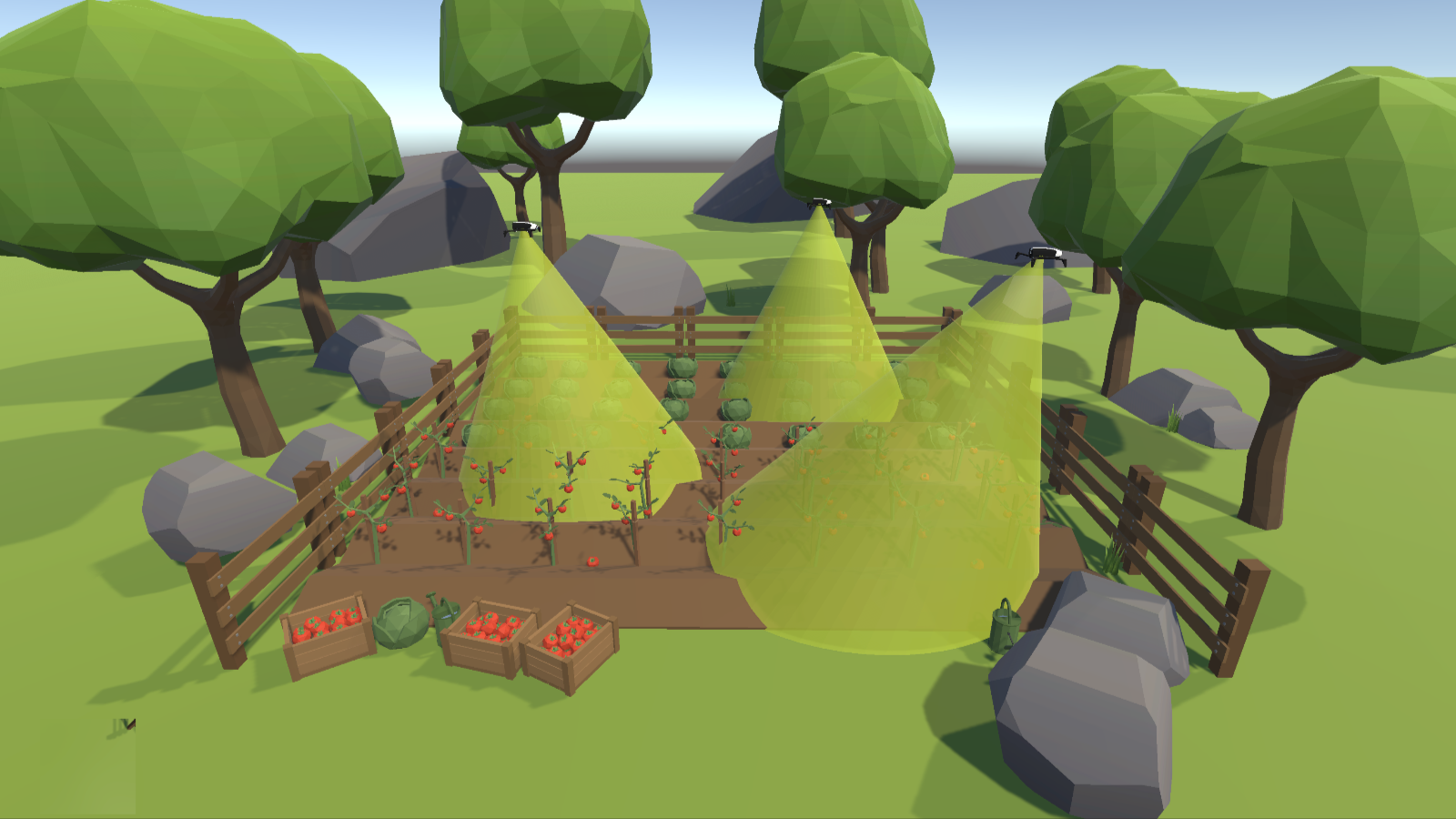}
    } \\
    \subfloat[Image captured by each drone]{
        \label{fig:cams}
        \resizebox*{4.5cm}{!}{\includegraphics{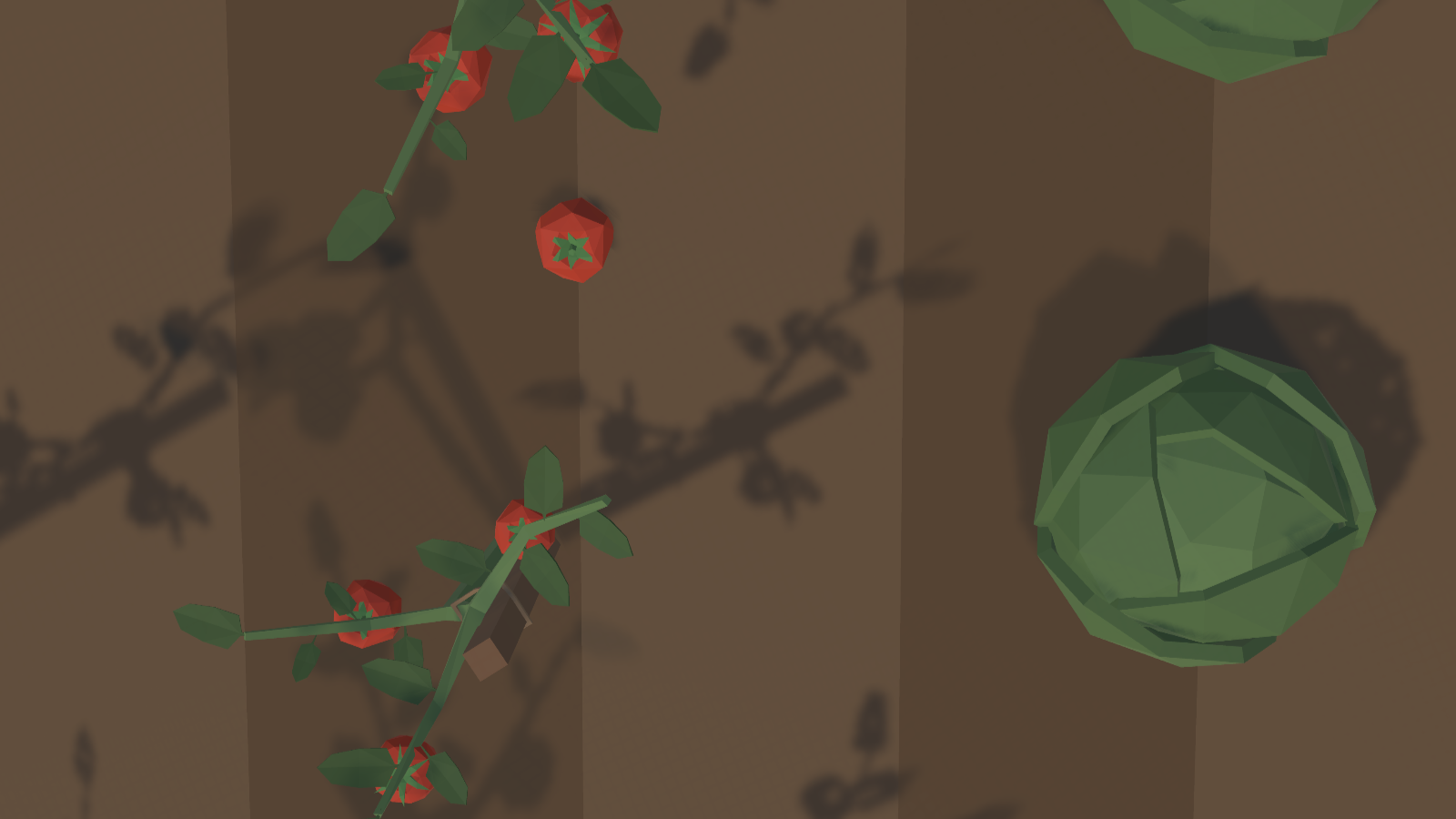}}\hspace{5pt}
        \resizebox*{4.5cm}{!}{\includegraphics{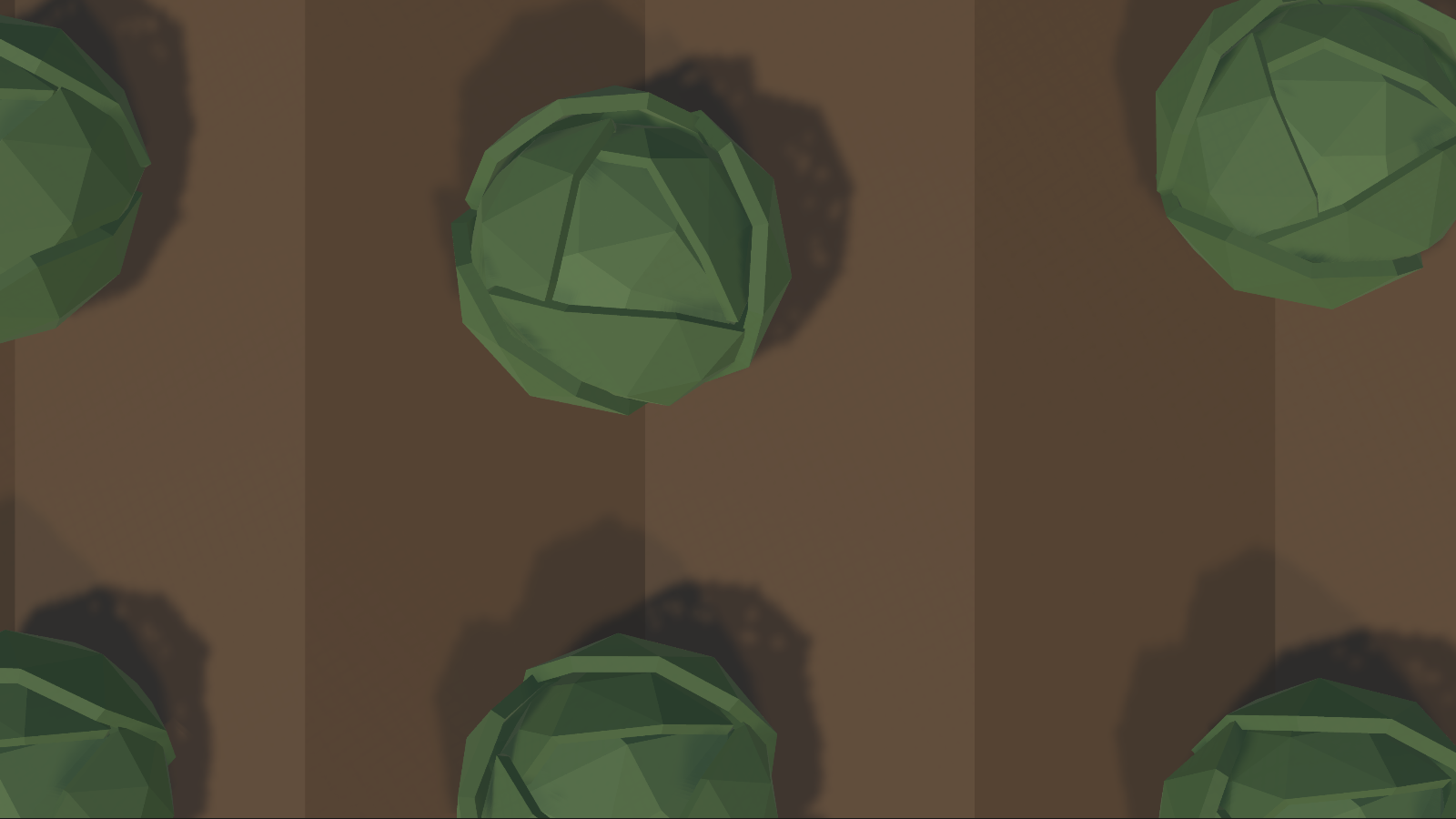}}\hspace{5pt}
        \resizebox*{4.5cm}{!}{\includegraphics{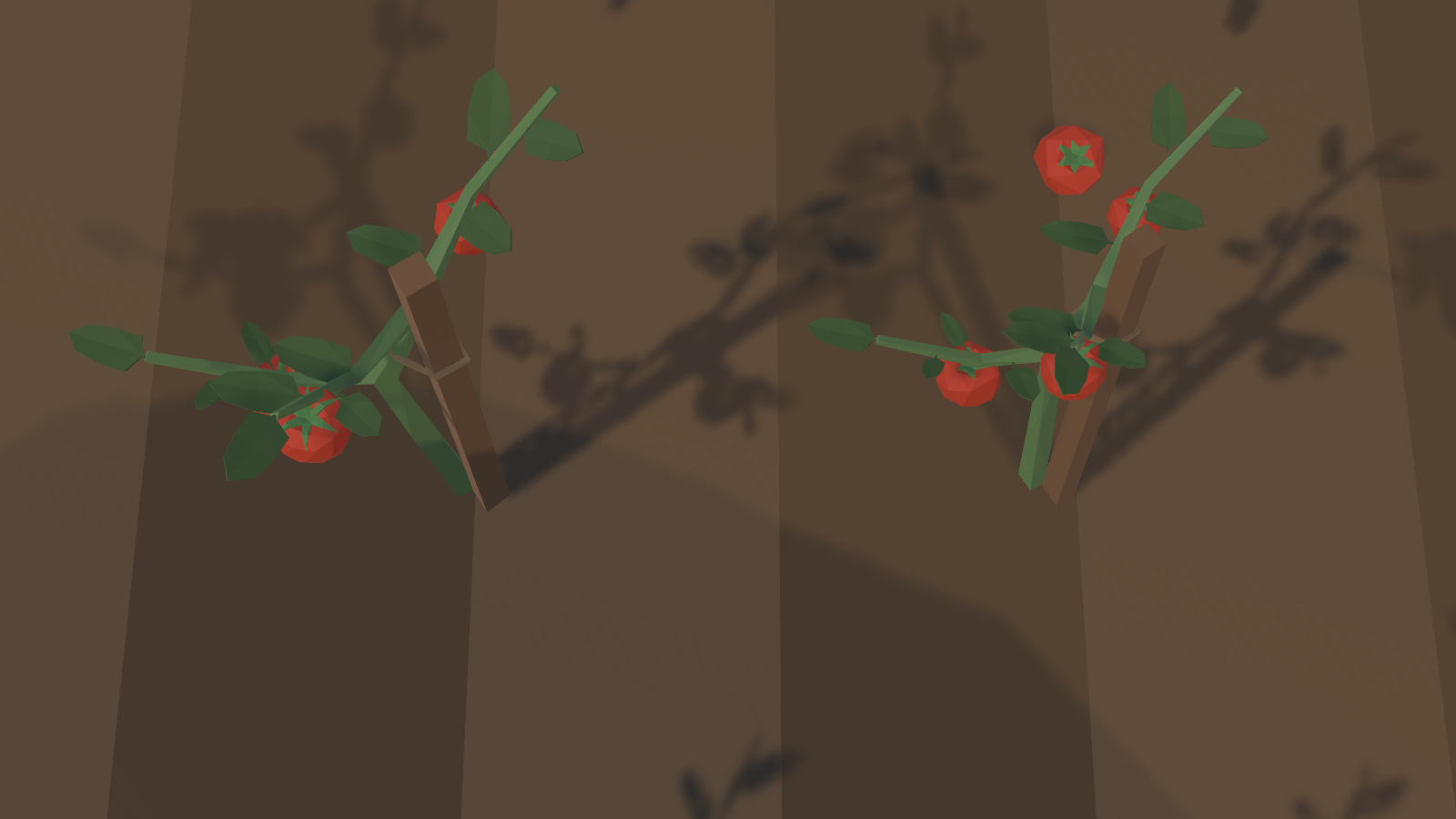}}\hspace{5pt}
    }
    \caption{Collaborative 3D map reconstruction using drone networks scene.} 
    \label{fig:scene}
\end{figure}

\section{Problem settings}

\subsection{Drones, Virtual Field, and Geometry}

\begin{figure}[t]
	\centering
	\includegraphics[scale=0.40]{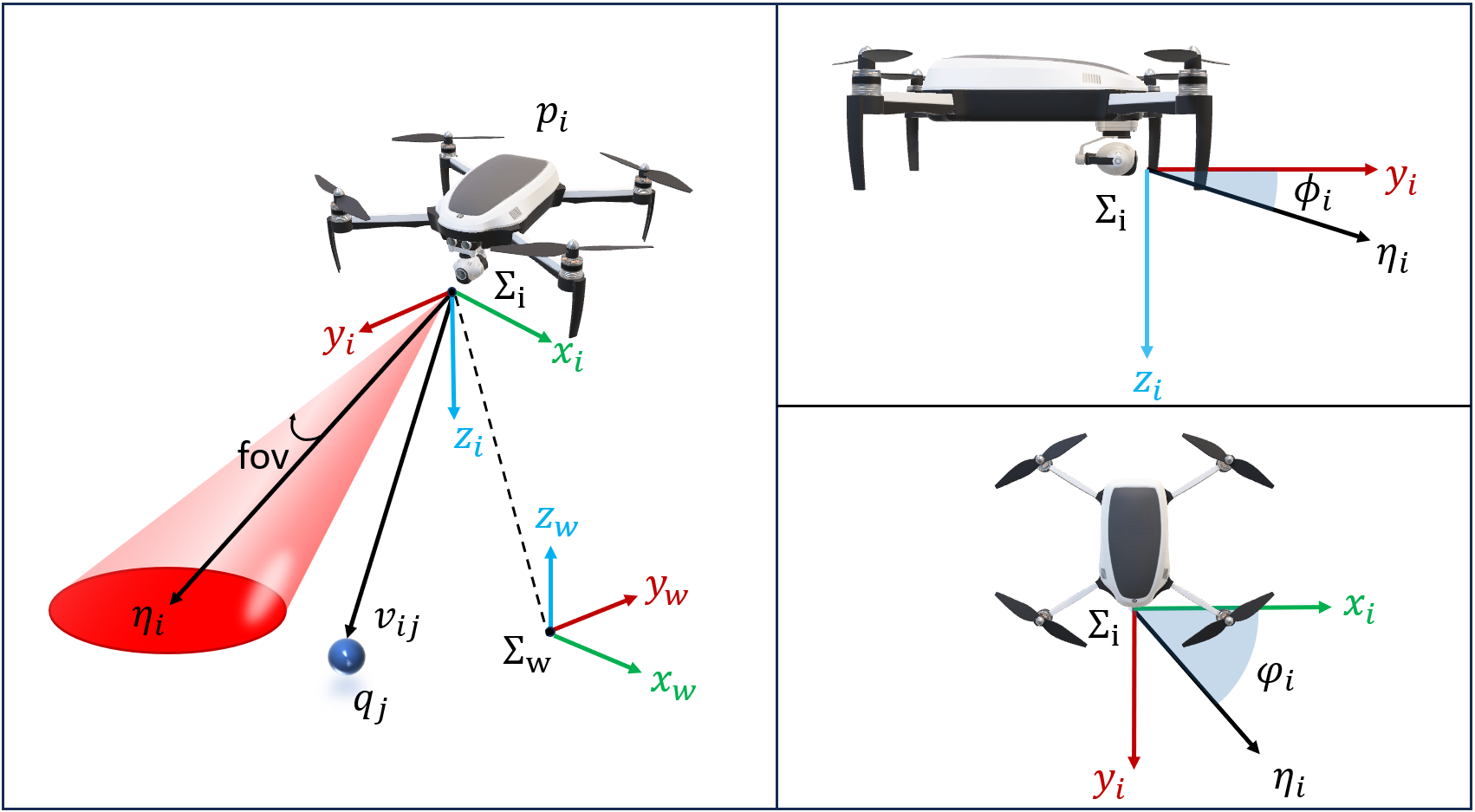}
	\caption{Gimbal coordinate reference.}
	\label{fig:coordinate_frame}
\end{figure}

We consider a scenario involving $n$ drones, each represented by the index set $\mathcal{I}:=\{1, \ldots, n\}$, operating in three-dimensional Euclidean space. All the drones are regarded as rigid bodies, and we define two frames of reference: the world frame $\Sigma_w$ and the body frame $\Sigma_i$, which is fixed on the camera of the $i$-th drone as illustrated in Fig~\ref{fig:coordinate_frame}. 
The $x$, $y$, and $z$ coordinates of the origin of $\Sigma_i$ relative to $\Sigma_w$ are denoted by $x_i$, $y_i$, and $z_i$, respectively.
The orientation of $\Sigma_i$ relative to $\Sigma_w$ is denoted by $R_i\in SO(3)$.
Each drone $i$ is equipped with an onboard camera that can be adjusted both horizontally, denoted by $\varphi_i$, and vertically, denoted by $\phi_i$, using a gimbal system. The horizontal angle $\varphi_i$ ranges from $0$ to $2\pi$, while the vertical angle $\phi_i$ ranges from $0$ to $\pi/2$, effectively forming a hemisphere of observable angles for the camera. 

During operation, all drones are assumed to maintain a constant altitude $z_i = z_c$. Consequently, we exclude the altitude component from the drone's position description, defining a state vector of drone $i$, $p_i := [x_i \ y_i \ \varphi_{i} \ \phi_{i}]^T \in \mathcal{P}\times [0, 2\pi) \times [0, \pi/2] \subset\mathbb{R}^4$, where $\mathcal{P}\subset \mathbb{R}^2$ represents a compact subset of a plane. 
Furthermore, each drone within the set $\mathcal{I}$ follows the subsequent dynamics: 

\begin{equation}
    \dot p_i = \left[\begin{array}{c}
        \dot x_i \\
        \dot y_i \\
        \dot \varphi_{i} \\
        \dot \phi_{i}
    \end{array}\right]
    =
    \left[\begin{array}{c}
        u^x_{i} \\
        u^y_{i} \\
        u^\varphi_{i} \\
        u^\phi_{i}
    \end{array}\right]
    = u_i.
\end{equation}
Here, $u^x_{i}$ and $u^y_{i}$ represent the linear velocity input for drone $i$, while $u^\varphi_{i}$ and $u^{\phi}_{i}$ correspond to the angular velocity input for adjusting the gimbal angles of drone $i$.

\begin{figure}[t]
	\centering
	\includegraphics[width=110mm]{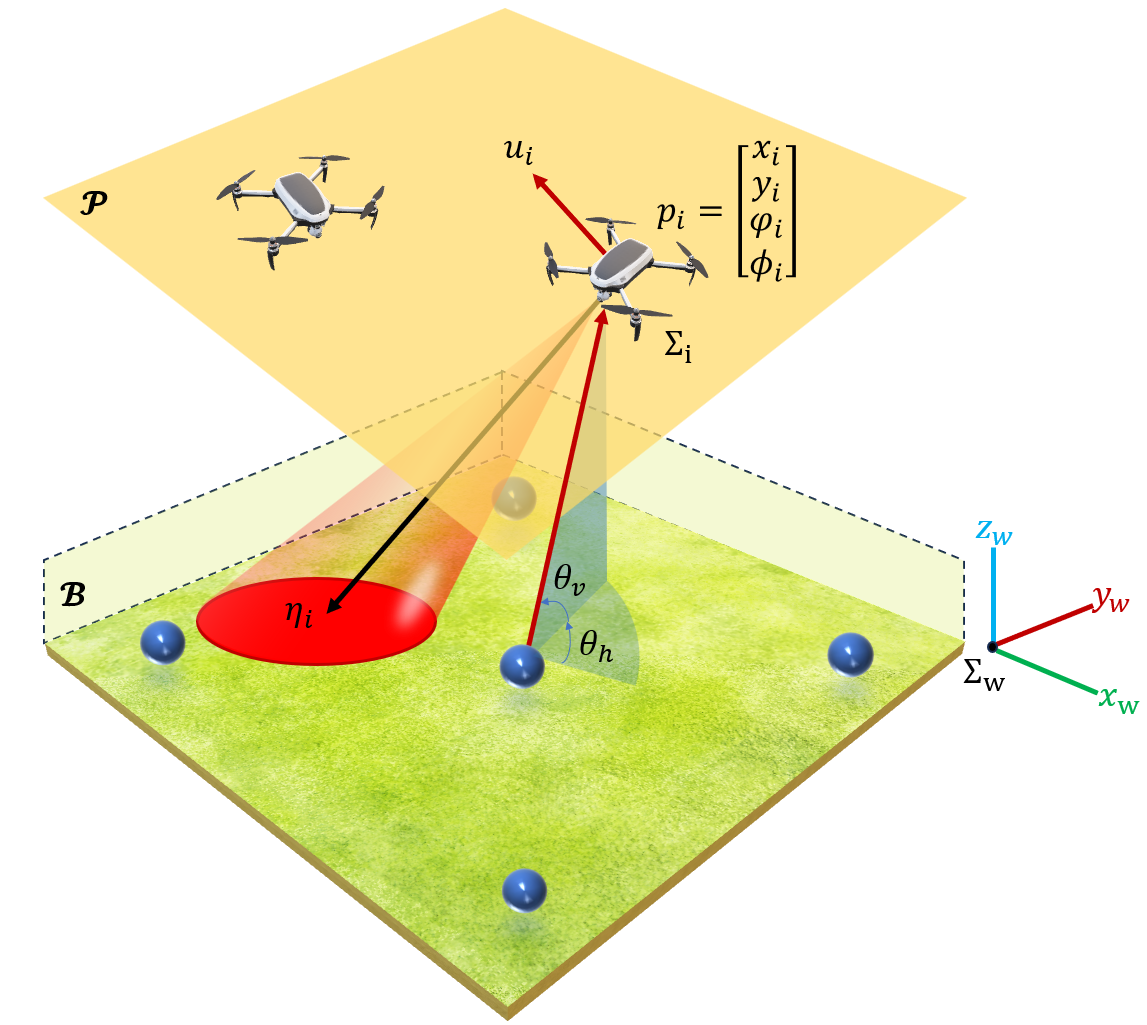}
	\caption{Illustration of the angle aware coverage problem with camera orientation control.}
	\label{fig:problem_illustration}
\end{figure}

Our main objective is to reconstruct 3D structure about a given target field. 
In this paper, we model the set of position coordinates in $\Sigma_w$ of all points in the target field to be observed as $\mathcal{B} \subset \mathbb{R}^3$, encompassing the highest and lowest points and any objects within. 
Now, it is widely known that every point in $\mathcal{B}$ should be observed by rich viewing angles in order to obtain a high-quality 3D structure.
To reflect this objective, we define the horizontal angle $\theta_h \in [-\pi, \pi)$ and the vertical angle $\theta_v \in (0, \pi/2]$ to represent the angles from which we observe specific target points. 
The points to be observed are then characterized by five variables consisting of not only the position coordinates $[x \ y \ z]^T \in \mathcal{B}$ but also the viewing angles $\theta_h$ and $\theta_v$ (See Fig~\ref{fig:problem_illustration}). Accordingly, we consider a coverage control problem over the 5D virtual space, $\mathcal{Q}_c \subset \mathbb{R}^5$, which encompasses all observation variables: $(x, y, z, \theta_h, \theta_v)$. Hence, the primary goal is how we can effectively control drone $p_i \in \mathcal{P}$ to monitor every point $q \in \mathcal{Q}_c$.
In the sequel, we discretize the 5D virtual space $\mathcal{Q}_c$ into $m$ cells and a representative point of the $j$-th cell is denoted by $q_j = [x_j\ y_j\ z_j\ \theta^h_j\ \theta^v_j] \in \mathcal{Q}_c$. 
The collection of $q_j, j=1,2,\dots,m$ is denoted by $\mathcal{Q}$.

Let us now introduce geometry associated with the present control problem.
Consider drone $i$ having the state $p_i := [x_i \ y_i \ \varphi_{i} \ \phi_{i}]^T$.
The optical axis of the camera is then described in $\Sigma_i$ as $\eta_i = [\cos(\varphi_i)\cdot\cos(\phi_i), \sin(\varphi_i)\cdot\cos(\phi_i), \sin(\phi_i)]^{\top}$.
On the other hand, the position vector in  $q_j$ (the first three elements) is described as 
\begin{equation}
v_{ij} = R_i^T v^w_{ij},\ v^w_{ij} = 
\begin{bmatrix}
x_j - x_i \\
y_j - y_i \\
z_j - z_c
\end{bmatrix}
\end{equation}
in $\Sigma_i$.
Also, the vector specified by the angles $\theta^h_{j}$ and $\theta^v_j$ is represented in $\Sigma_w$ as 
\begin{equation}
v_{j} = \begin{bmatrix}
\cos(\theta_h)\cdot\cos(\theta_v)\\ \sin(\theta_h)\cdot\cos(\theta_v)\\
\sin(\theta_v)
\end{bmatrix} .
\end{equation}.

\begin{remark}
In the context of drone camera capabilities, it is essential to address certain limitations that may affect their performance. While some drone cameras offer pitch, yaw, and roll rotation capability, the structural limitations of many gimbals often restrict the range of yaw rotation angles, often preventing a complete 360-degree rotation. To address this, controlling the camera's yaw angle typically involves coordination with the drone's body rotation. However, in this paper, we do not focus on the drone's rotational behavior, and for simplicity, we consider it as part of the gimbal's angle control. Moreover, it's important to note that camera sensors generally have rectangular designs, leading to a rectangular beam-shaped observable area for the drone. Nevertheless, for the sake of simplicity and without loss of generality, we assume a circular shape for the drone's observable area.
\end{remark}

\subsection{Performance function and importance index}

In this subsection, we present a performance function, denoted as $h: \mathcal{P}\times [0, 2\pi) \times [0, \pi/2] \times \mathcal{Q}_c \rightarrow [0, 1]$. This function characterizes a drone's coverage capability at state $p_i$ concerning a specific point $q_j \in \mathcal{Q}$. A higher value of the performance function corresponds to a more effective coverage capability. The effective coverage of point $q_j$ by a drone positioned at state $p_i$ depends on two key factors:

\begin{enumerate}
    \item The geometric position of point $q_j$ must fall within the drone's camera's field of view.
    \item The relative geometric relationship between the drone's position and the observed point $q_j$ must satisfy the angle $(\theta_h,\theta_v)$ of $q_j$.
\end{enumerate}

The first condition can then be described as follows:

\begin{equation}
   h_1(p_i,q_j) = (\text{fov} - \arccos(\eta_i \cdot v_{ij}/\|v_{ij}\|))^2 \geq 0,
\end{equation}
where $\text{fov}$ is the half-angle of the camera's field of view.
The second condition can be described in a similar way:
\begin{equation}
   h_2(p_i,q_j) = (\arccos(v_j \cdot v^w_{ij}/\|v_{ij}^w\|))^2 \geq 0. 
\end{equation}

In the subsequent controller design, we employ the following performance function combining $h_1$ and $h_2$ with appropriate tuning of the parameters $\sigma_1$ and $\sigma_2$:
\begin{equation}
    h\left(p_i, q_j\right):=\exp \left(-\frac{h_1(p_i,q_j)}{2 \sigma_1^2}\right) \exp \left(-\frac{h_2(p_i,q_j)}{2 \sigma_2^2}\right) .
\end{equation}

Remark that the parameter $\sigma_1$ should be tuned so that the performance function is sufficiently close to 0 for any point $q$ outside of the field of view. On the other hand, we have not thought of any systematic way to tune $\sigma_2$, and tune empirically this parameter in the subsequent simulation. We would like to leave the issue to future work.

Let us next introduce the importance index $\psi_j\in [0, \infty)$ assigned to each point $q_j \in \mathcal{Q}$. 
From the definition of the performance function $h$, a large value of $h(p_i,q_j)$ for some $i$ means that
the drone already samples a good image on $q_j$ and then we can reduce the importance of the point $q_j$.
Based on this observation, similarly to \cite{shimizu2021angle},  we propose the following update law for the importance index $\psi_j$:

\begin{equation}
    \dot{\psi}_j=-\delta \max _{i \in \mathcal{I}} h\left(p_i, q_j\right) \psi_j\ \ (\delta>0) .
    \label{eqn:update_psi}
\end{equation}

In terms of the controller design, it is preferable that the performance function meets these two properties:
\begin{itemize}
    \item Restricting the performance function within the range from 0 to 1 in order to govern the decay rate of $\psi_j$ in (7). 
    \item Ensuring a continuous non-zero gradient for the performance function in order to ensure the existence of the input $u_i$ that increases the function $b_i$.
\end{itemize}
An example that meets the above properties is the Gaussian function. This is why we utilize it as the performance function.

\subsection{Objective function}

We are now prepared to present the objective function that needs to be minimized, defined as

\begin{equation} \label{eq:J}
    J:= \sum_{j=1}^m \psi_j .
\end{equation}
This equation represents the integral of the density function across the entire region. When the density function approaches values close to $0$ for individual points, it signifies the effective capture of images for those points. As $J$ approaches $0$, it indicates that the region's coverage is nearing completion, which is crucial for capturing images suitable for reconstructing the 3D structure. Consequently, the primary objective is to steer the drones in a manner that drives $J$ towards convergence to zero. 
However, due to the monotonically decreasing property of $\psi_j$ in (\ref{eqn:update_psi}), it is trivial to achieve $J \to 0$ itself.
We thus impose another specification 
$\dot J \leq -\gamma$ for a positive constant $\gamma$ to specify efficiency of the task completion. This covers the primary objective $J \to 0$.
In summary, the problem to be addressed in this paper is to determine the velocity input $u_i$ so that $\dot J \leq -\gamma$ is satisfied for a given parameter $\gamma > 0$.

\section{QP-based Controller Design}

In this section, we propose a controller based on quadratic programming, which enforces the constraint $\dot J \leq -\gamma$ using the concept of control barrier functions. To this end, we begin by computing the time derivative of $J$:

\begin{align}
\dot{J} &= \sum_{j=1}^m \dot \psi_j(q) = - \sum_{j=1}^m \delta \max _{i \in \mathcal{I}} h\left(p_i, q_j\right) \psi_j  \nonumber\\
    &= -\sum^{n}_{i=1} \int_{j\in \mathcal{V}_{i}(p)} \delta h(p_{i},q_{j}) \psi_j = - \sum^{n}_{i=1} I_{i} ,
\end{align}
where 
\begin{align} \label{eq:I_i}
	I_{i} := \int_{j\in \mathcal{V}_{i}(p)} \delta h(p_{i},q_{j})\psi_{j}
\end{align}
corresponds to the contribution by drone $i$ to reduce $J$ in (\ref{eq:J}). The set $\mathcal{V}_{i}(p)$, which depends on $p := (p_i)_{i\in \mathcal{I}}$, is a Voronoi-like partition of the set $\mathcal{M} := \{1, 2, \dots, m\}$ defined as 
\begin{align}
	\mathcal{V}_{i}(p) := \{j \in \mathcal{M} \ | h({p_{i}, q_{j}}) \leq h({p_{k}, q_{j}}) \ \forall k \in \mathcal{I}\}.
\end{align}

By defining $b_{i,I} = I_i - \gamma/n$ as a candidate for the control barrier function, with a given $\gamma > 0$, we aim to compel the drone to minimize the objective function $J$ at a specified rate of decrease, denoted by $\gamma$. This criterion is met when the control input $u_i$ is chosen such that the following inequality holds:

\begin{align}
	\left( \frac{\partial b_{i,I}}{\partial \psi_j} \right) \dot \psi_j + \left( \frac{\partial b_{i,I}}{\partial p_i} \right)^T u_i + \alpha_1(b_{i,I}) \geq 0 ,
\end{align}
where $\alpha_1$ is a locally Lipschitz extended class $\mathcal{K}$ function. A continuous function
$\alpha : (- b, a) \xrightarrow[]{} (-\infty, \infty)$ 
is said to belong to extended class $\mathcal{K}$ for some $a, b > 0$ if it is strictly increasing and $\alpha(0) = 0$ \cite{ames2017control}.

Furthermore, to adhere to the gimbal input limitations, we introduced a secondary control barrier function constraint, restricting the vertical angle $\phi_i$ of the drone's gimbal within the range of $0$ to $\pi/2$. While functionally equivalent vertical angles exist between $\pi$ and $\pi/2$ as those between $0$ and $\pi/2$, the physical structure of the gimbal imposes constraints on this range. Failing to enforce these constraints during optimization may result in commands that surpass the $\pi/2$ threshold, which not only prevents gimbal from executing the command, but also results in a reduction in horizontal rotation. By maintaining the vertical angle within the specified limits, we ensure optimal performance and prevent unintended consequences associated with exceeding the gimbal's mechanical boundaries.

Therefore, to ensure the vertical angle $\phi_i$ of the gimbal varies between $\phi_{\text{{min}}}$ and $\phi_{\text{{max}}}$, we define a control barrier function. This function, denoted as $b_{i,\phi}$, ensures that the condition

\begin{equation}
    b_{i,\phi}= \left(\phi_{\text{{max}}} - \phi_{\text{{min}}}\right)^2 - \left( \phi_{i} - \frac{\phi_{\text{{min}}} + \phi_{\text{{max}}}}{2} \right)^2 \geq 0
\end{equation}
is satisfied at all times, securing the appropriate vertical angle constraints throughout the optimization process. This criterion is met when the control input $u_i$ is chosen such that the following inequality holds:

\begin{align}
	 \left( \frac{\partial b_{i,\phi}}{\partial p_i} \right)^T u_i + \alpha_2(b_{i,\phi}) \geq 0 ,
\end{align}
where $\alpha_2$ is a locally Lipschitz extended class $\mathcal{K}$ function.

Now, we are ready to present the QP-based controller as follows:

\begin{subequations} \label{eq:qp_controller}
	\begin{align}
		(u_{i}^{*},w_{i}^{*}) &= \argmin_{ (u_{i},w_{i}) \in \mathcal{U} \times \mathbb{R}} \epsilon\norm{u_{i}}^2 + |w_{i}|^2  \label{eq:qp_controller_1} \\
		&\text{s.t.} \quad \dot b_{i,I}+\alpha_1\left(b_{i,I}\right) \geq w_i ,\label{eq:qp_controller_2} \\
            & \quad \quad \ \dot b_{i,\phi} +\alpha_2\left(b_{i,\phi}\right) \geq 0 ,\label{eq:qp_controller_3}
	\end{align}
\end{subequations}
where $\epsilon$ is the penalty variable and  $w_i$ is the slack variable. 

\begin{theorem} \label{th:qp_form}
    Suppose that no $q_j (j \in \mathcal{M})$ is located on the boundary of $\mathcal{V}_i(p)$. When $\alpha_1: \mathbb{R} \rightarrow \mathbb{R}$ and $\alpha_2: \mathbb{R} \rightarrow \mathbb{R}$ are set as a linear function such that $\alpha_1(b_{i,I})=a_1 b_{i,I}$ and $\alpha_2(b_{i,\phi})=a_2 b_{i,\phi}$, where $a_1>0$ and $a_2>0$ is a positive scalar, the problem is equivalently reformulated as :

    \begin{subequations} \label{eq:qp_controller_final}
	\begin{align}
		(u_{i}^{*},w_{i}^{*}) &= \argmin_{ (u_{i},w_{i}) \in \mathcal{U} \times \mathbb{R}} \epsilon\norm{u_{i}}^2 + |w_{i}|^2  \\
		&\text{s.t.} \quad \xi_{1i}^{T}u_{i} + \xi_{2i} \geq w_{i} , \\
            & \quad \quad \ \chi_{1i}^{T}u_{i} + \chi_{2i} \geq 0 ,
	\end{align}
    \end{subequations}
    where

    \begin{subequations}
	\begin{align}
		& \xi_{1 i}:=\int_{j \in \mathcal{V}_i(p)}\delta \frac{\partial  h\left(p_i, q_j\right)}{\partial p_i} h\left(p_i, q_j\right) \psi_j , \\
            & \xi_{2 i}:=- \frac{a_1 \gamma}{n}+\int_{j \in \mathcal{V}_i(p)}\left(-\delta^2 h^2\left(p_i, q_j\right)+a \delta h\left(p_i, q_j\right)\right) \psi_j ,
	\end{align}
    \end{subequations}

    \begin{subequations}
	\begin{align}
		& \chi_{1 i}:= \left[\begin{array}{c}
        0 \\
        0 \\
        0 \\
        2 \left ( \phi_{i} - \frac{\phi_{\textup{{min}}} + \phi_{\textup{{max}}}}{2} \right)
    \end{array}\right] , \\
            & \chi_{2 i}:= 
             a_2 \left ( \left(\phi_{\textup{{max}}} - \phi_{\textup{{min}}}\right)^2 - \left( \phi_{i} - \frac{\phi_{\textup{{min}}} + \phi_{\textup{{max}}}}{2} \right)^2 \right ) .
	\end{align}
    \end{subequations}

\end{theorem}
Please refer to Appendix \ref{app: theorem_1} for the detail proof of the controller (\ref{eq:qp_controller_final}). 

\begin{remark}

The optimization problem (15) is always feasible because the constraint is softened by the 
 slack variable $w_i$. However, the existence of the solution to (15) does not mean that the original specification $I_i \geq \gamma/n$ is met.
For a too large $\gamma$, the specification might not be satisfied in the transient in the presence of the velocity limit of the drone.
It is desirable that $\gamma$ is appropriately determined so that $I_i \geq \gamma/n$ is violated only when the coverage is almost completed, namely $J \approx 0$. However, a systematic way to determine such $\gamma$ for a given velocity limit is left for future work.

\end{remark}

\section{Computation acceleration}

In contrast to conventional two-dimensional coverage control, dealing with coverage in a five-dimensional space significantly increases the computational load. Specifically, the number of cells, $m$, drastically increases in the five-dimensional case, which affects various processes including calculation of the performance function, updates of the importance index, Voronoi-like region partitioning, and, more significantly, gradient computation of the performance function. In particular, the gradient computation in terms of the five variables $x_i$, $y_i$, $\varphi_{i}$ and $\phi_{i}$ takes up most of the computational time.
\cite{shimizu2021angle} presented a computationally efficient implementation of the controller with slight approximations in computation,
where the position of the drone monitoring a point from a specified view angle was uniquely determined, enabling the mapping of the five-dimensional target field to a two-dimensional space with slight approximation. On the other hand, in this paper, we assume that the camera orientation can be controlled. In this case, the camera state monitoring a point from a specified view angle is not uniquely determined, and the mapping from five-dimensional field to a lower-dimensional space cannot be defined. Therefore, the same approach could not be applied to the present problem.
We thus approach this problem technologically, namely we introduce a JAX library accelerating the processing through just-in-time (JIT) compilation and GPU acceleration.

\subsection{JIT compilation}

As same as in \cite{shimizu2021angle}, our program is primarily written in Python, and we use NumPy library for matrix calculations. Since Python is an interpreted language, the code is not compiled before execution. This leads to a significant performance loss when running large-scale computations in Python compared to statically compiled code. The complex calculations involved in computing the gradient of the performance function exacerbate this loss, making real-time calculations challenging.
We thus replace NumPy by JAX's built-in NumPy, allowing to compile functions of pure matrix computations. 

\begin{figure}
    \centering
    \includegraphics[width=12cm]{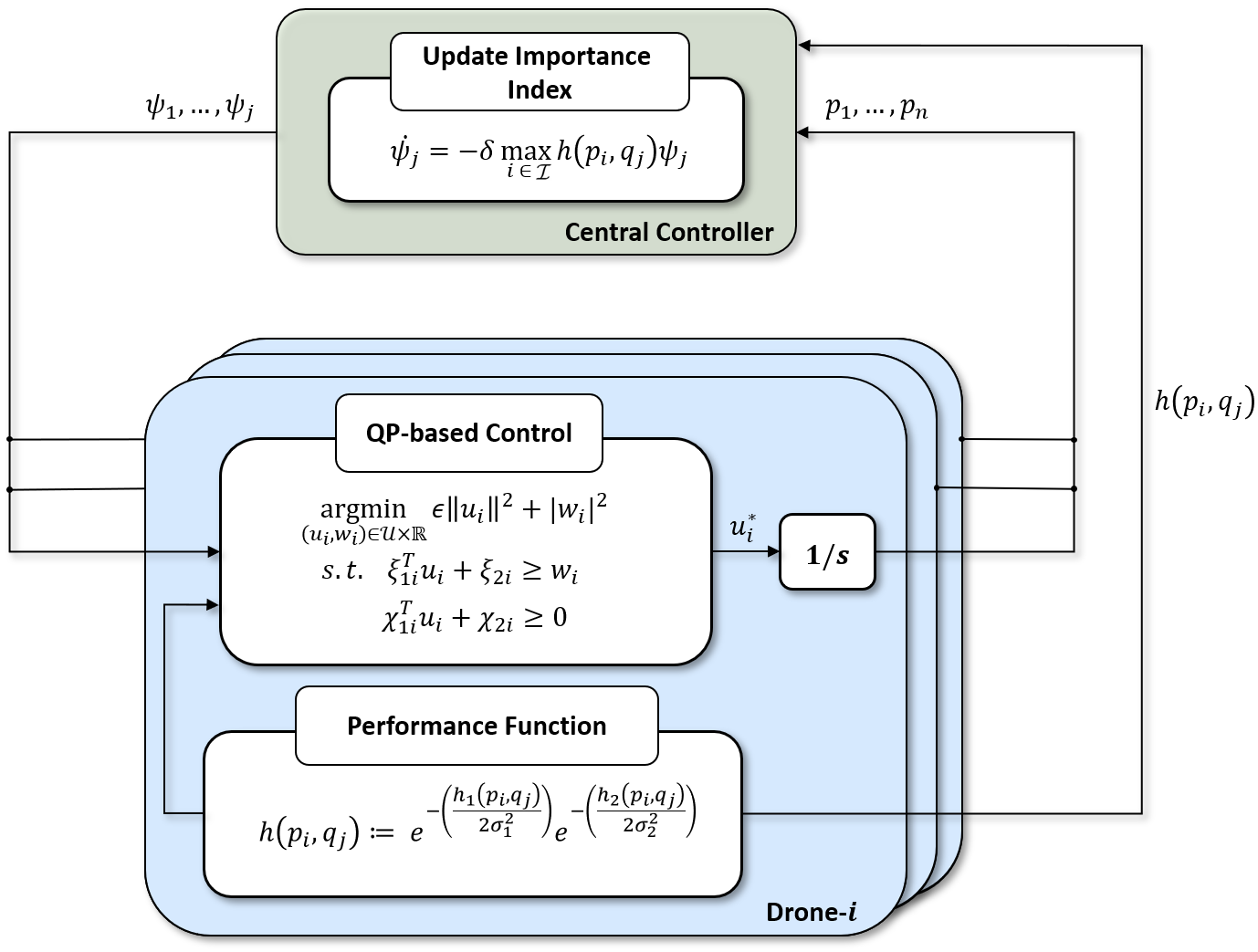}
    \caption{Controller architecture}
    \label{fig:controller_arch}
\end{figure}

We use JIT in two types of ROS nodes in our program: central controller and drone controller, where the importance index is updated in the central controller, while the computation performance function, dividing the Voronoi-like partition and the gradient computation is done in the drone controller, as Fig. \ref{fig:controller_arch} shows. The first two rows of table \ref{table:computation_time} compare the impact on the computation speed of whether JIT is used or not. For drone controller, which is more computationally intensive compared to central controllers, computation times over 1s are unacceptable for real-time control.
Meanwhile, after accelerating the function using JIT, the computational speed gains a significant increase even with the same computer, enabling it to reach a computational frequency of around 10Hz. 

\begin{table}[h!]
\setlength{\tabcolsep}{1.6mm}
\tbl{Average computation time (ms) of one step}
{
\begin{tabular}{ |c|c|c| }
\hline
  & Central controller & Drone controller\\
\hline
CPU & 212 & 1206  \\ 
\hline
JIT, CPU & 31 & 102 \\
\hline
JIT, GPU & 6 & 22 \\
\hline
\end{tabular}
}
\label{table:computation_time}
\end{table}

\subsection{GPU acceleration}

\begin{figure}
\centering
    \subfloat[drone controller]{\resizebox*{7cm}{!}{\includegraphics{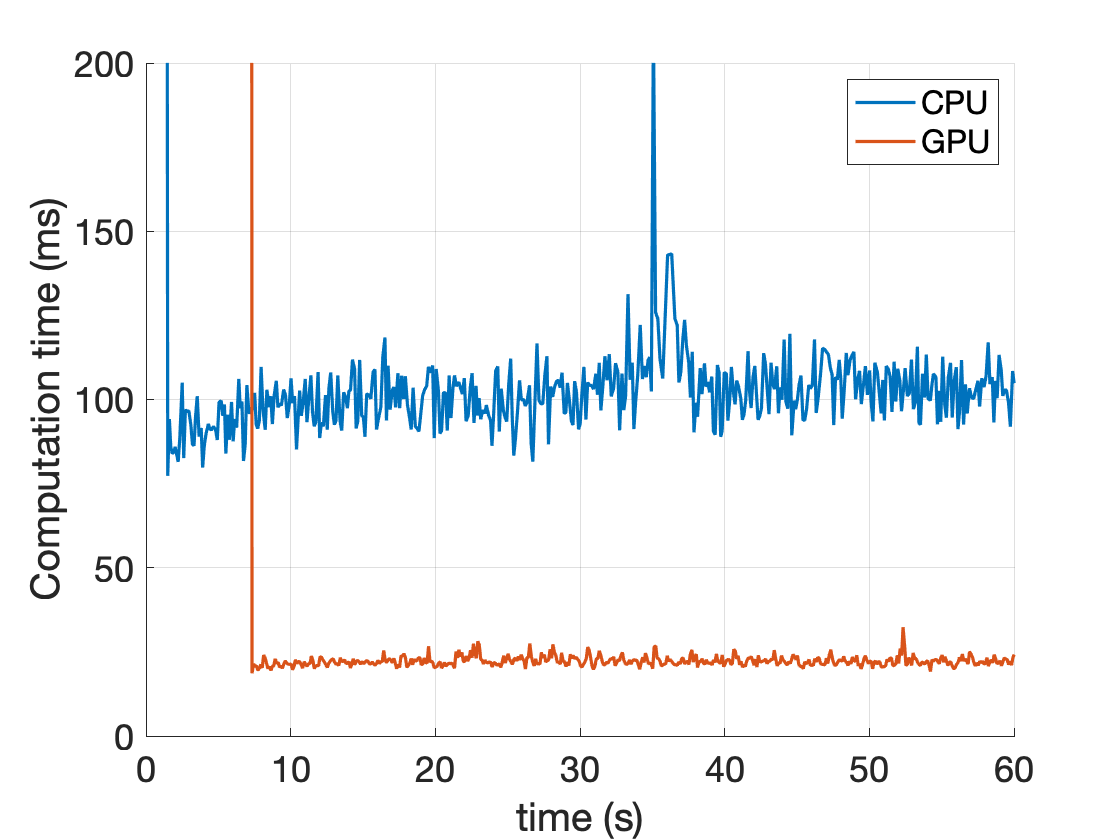}}}\hspace{5pt}
    \subfloat[central controller]{\resizebox*{7cm}{!}{\includegraphics{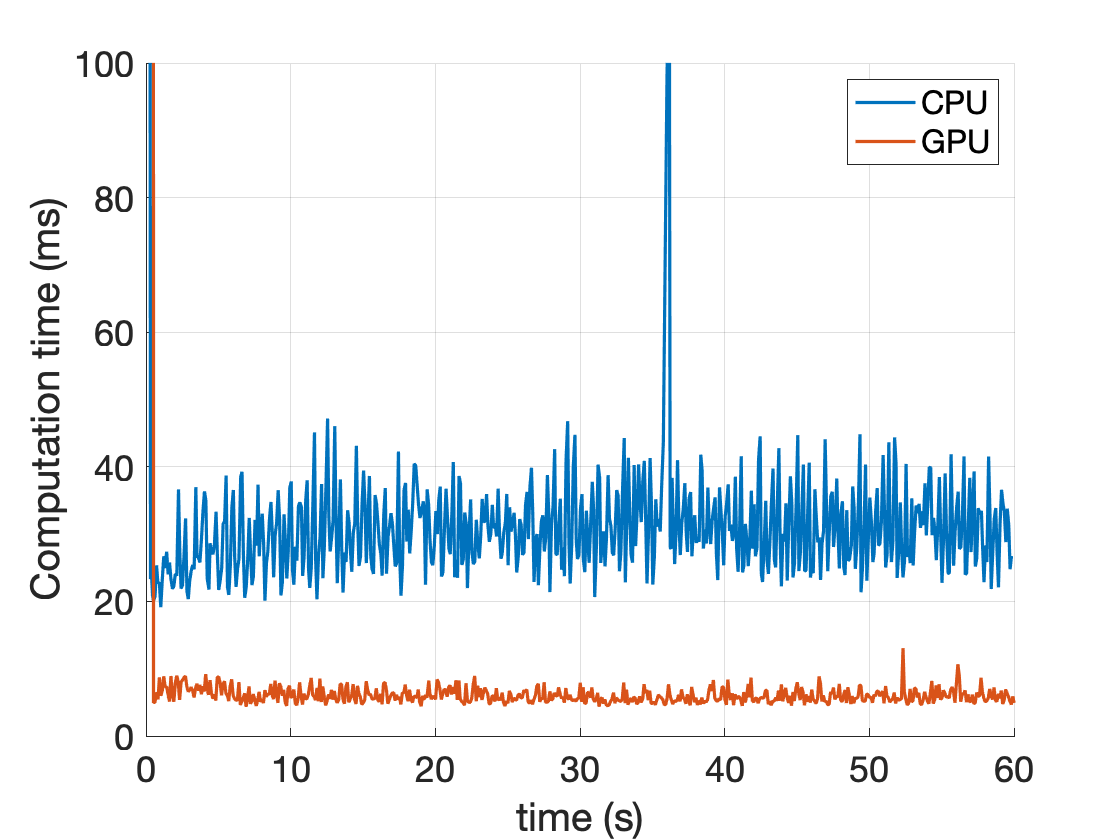}}}\hspace{5pt}
    \caption{Comparisons of CPU and GPU computational time.}
    \label{fig:benchmark}
\end{figure}

Since CPU performance is very difficult to improve and the computation time using CPU can be greatly altered by the size of the matrix, we utilize JAX's GPU acceleration feature to do matrix calculations on GPU.

We conducted a timing comparison on a laptop running Ubuntu 20.04.3 LTS, comparing the time it takes to perform computations using the CPU and GPU. The laptop is equipped with a 20-core Intel i7-12800H CPU and a Nvidia RTX A3000 GPU. Fig \ref{fig:benchmark} shows how the computation time consumed at each step changes over time when using CPU and GPU. We note that the first step of the computation takes a very long time in comparison, which is caused by the JIT compilation during the initial execution of functions. 
When using the GPU, compilation took longer, but after compilation was complete, computation using the GPU took roughly $1/5$ the time of the CPU, and the time spent was relatively stable in comparison, while the time spent on the CPU fluctuated considerably.

During CPU-based computations, the CPU usage is around 60\%, and during GPU-based computations, the GPU usage is reported only 24\%. As indicated in the table \ref{table:computation_time}, employing GPU acceleration for matrix computations results in significant time savings, particularly for large matrices, despite the overhead of data transfer between the CPU and GPU memory. It is to be expected that after increasing the size of the matrix, the computation time of the CPU and GPU will produce a bigger gap.

\section{Controller Evaluation and Verification}

\begin{table}[h!]
    \centering
    \setlength{\tabcolsep}{1.6mm}
    \caption{Parameter setting}
    \label{table:parameter}
    \begin{tabular}{ c c c c c c c }
        \toprule
        $a_1$ & $a_2$ & $\sigma_1$ & $\sigma_2$ & $\epsilon$ & $\gamma$ & $\delta$ \\
        \midrule
        5.0 & 1.0 & 0.13 & 0.18 & 0.0001 & 0.05 & 5.0  \\ 
        \bottomrule
    \end{tabular}
\end{table}

In this section, we demonstrate the proposed controller through simulation in the ROS Noetic environment.

The monitored space was set as a cube shape with the range of $[-1,1]\mathrm{m} \times [-1,1]\mathrm{m} \times [0,0.5]\mathrm{m}$. The viewing angle space is set to $\theta_h \in [-\pi, \pi)$, and $\theta_v \in [ \pi/6, \pi/2 ]$. We set the number of drones $n$ to 3 and their initial positions are uniformly distributed to $p_1 = [1.0, 0.2]^T \mathrm{m}$, $p_2 = [-1.0, -0.2]^T \mathrm{m}$ and $p_3 = [0.0, 0.5]^T \mathrm{m}$. Their flight altitude is fixed to $h_c = 1.0 \mathrm{m}$ after takeoff. To simulate a drone equipped with a camera with a focal length of $x$[mm], we set the field of view of the drone to $\pi/6$. The initial angles of the gimbal for the three drones are all set to $\varphi_i = 0$ and $\phi_i = \pi/2$, indicating that the cameras were initially pointed vertically downwards. The horizontal angles of the drone gimbal $\varphi_i$ were limited to $[0, 2\pi)$, while the vertical angles $\phi_i$ were limited to $(0, \pi/2]$.
The field is divided into $m = 1.5 \times 10^7$ small cells of size $0.02\mathrm{m} \times 0.02\mathrm{m} \times 0.1\mathrm{m} \times \pi/30 \mathrm{rad} \times \pi/30 \mathrm{rad}$. 
Other parameters are presented as table \ref{table:parameter}.
The quadratic programming is solved by CVXOPT \cite{andersen2013cvxopt}.

\begin{figure}[t]
\centering
    \subfloat[$t=0s$]{\resizebox*{4cm}{!}{\includegraphics{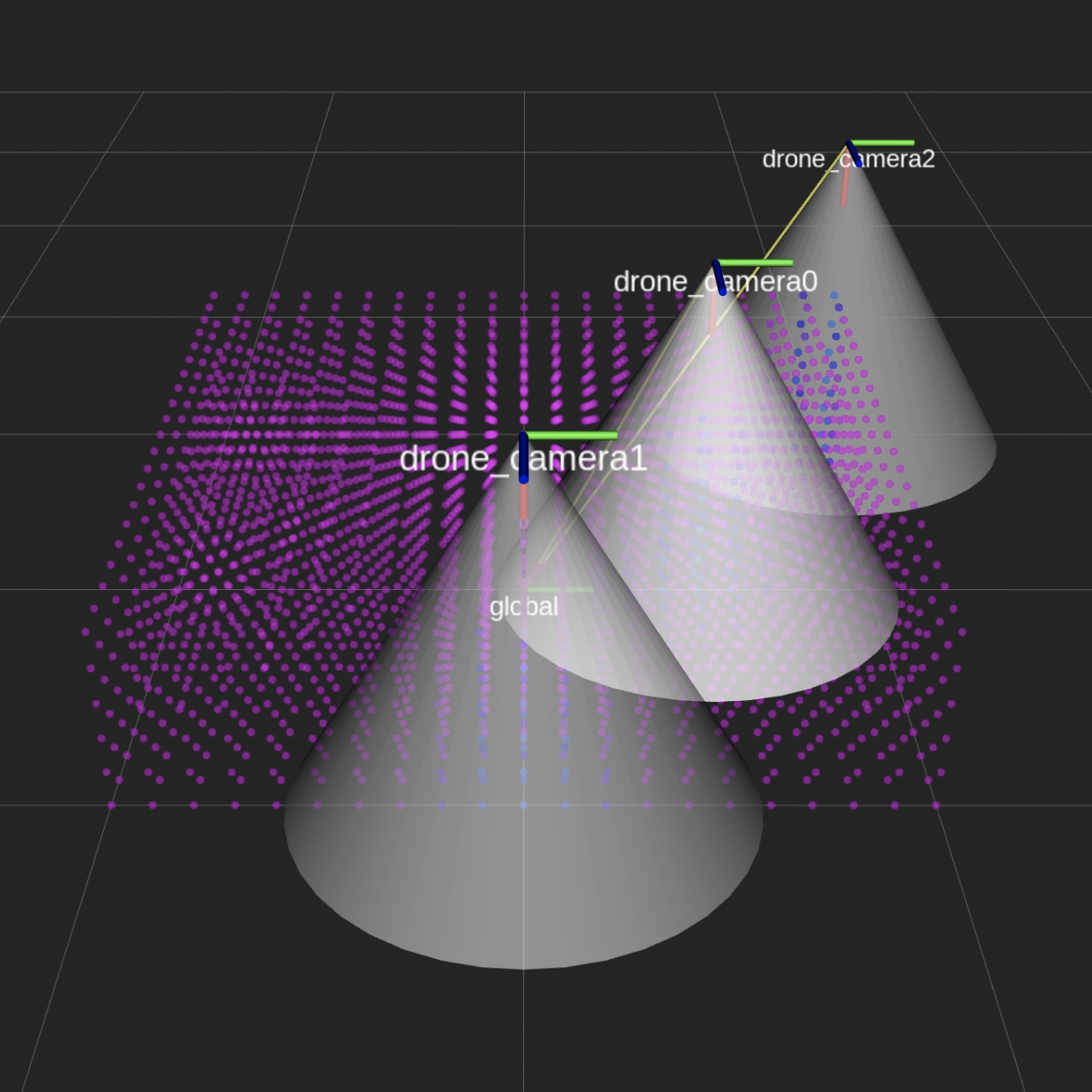}}}\hspace{5pt}
    \subfloat[$t=30s$]{\resizebox*{4cm}{!}{\includegraphics{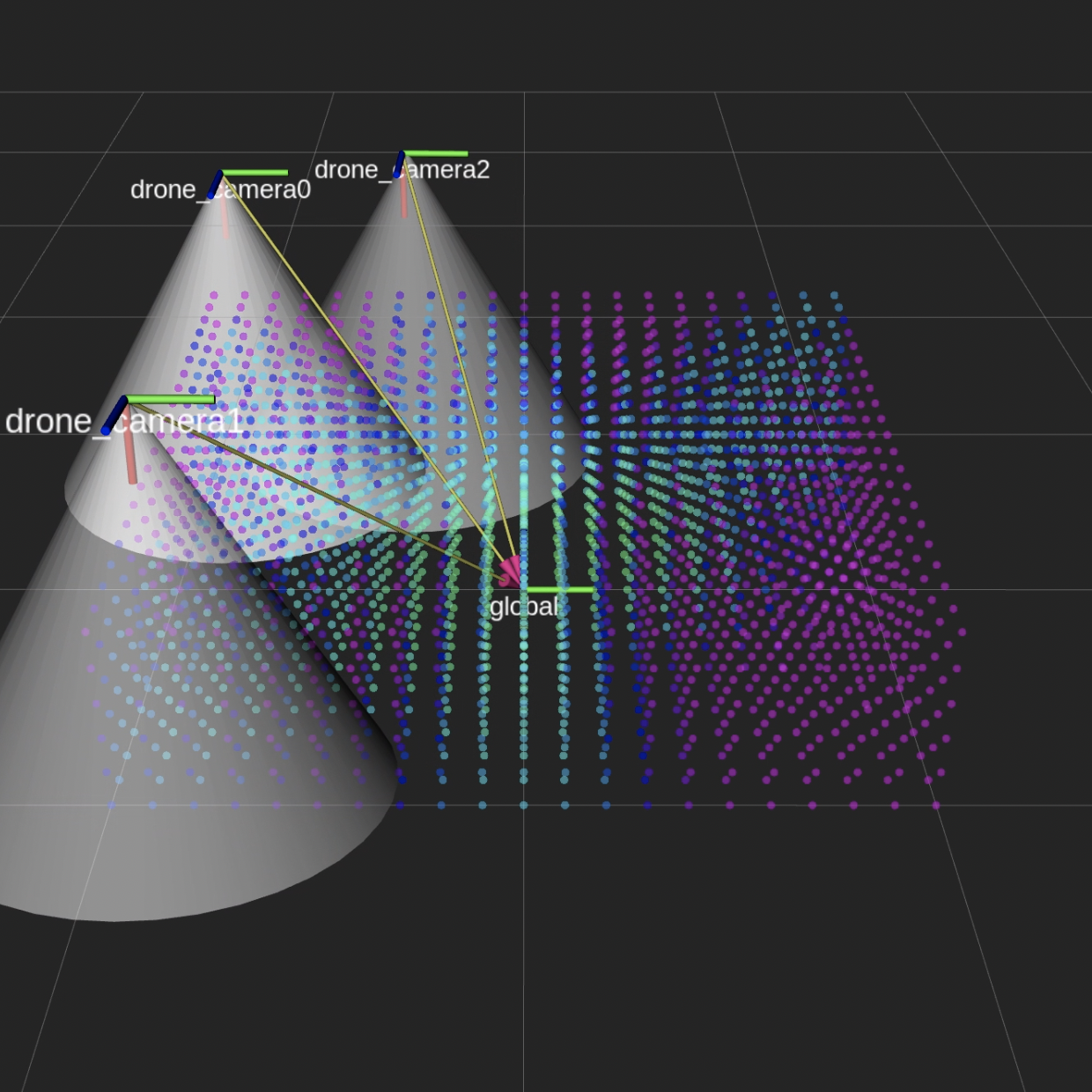}}}\hspace{5pt}
    \subfloat[$t=50s$]{\resizebox*{4cm}{!}{\includegraphics{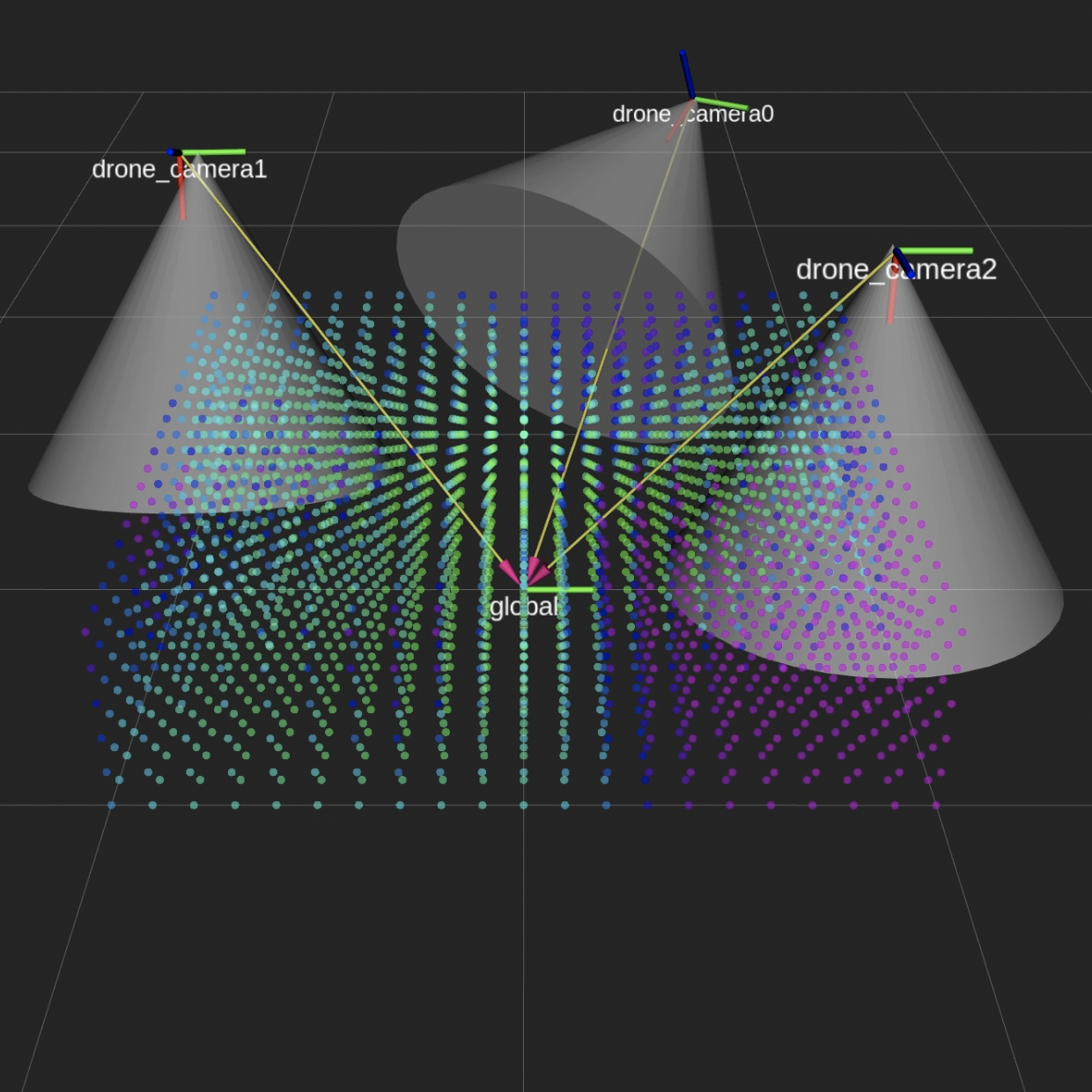}}}\hspace{5pt}
    \subfloat[$t=100s$]{\resizebox*{4cm}{!}{\includegraphics{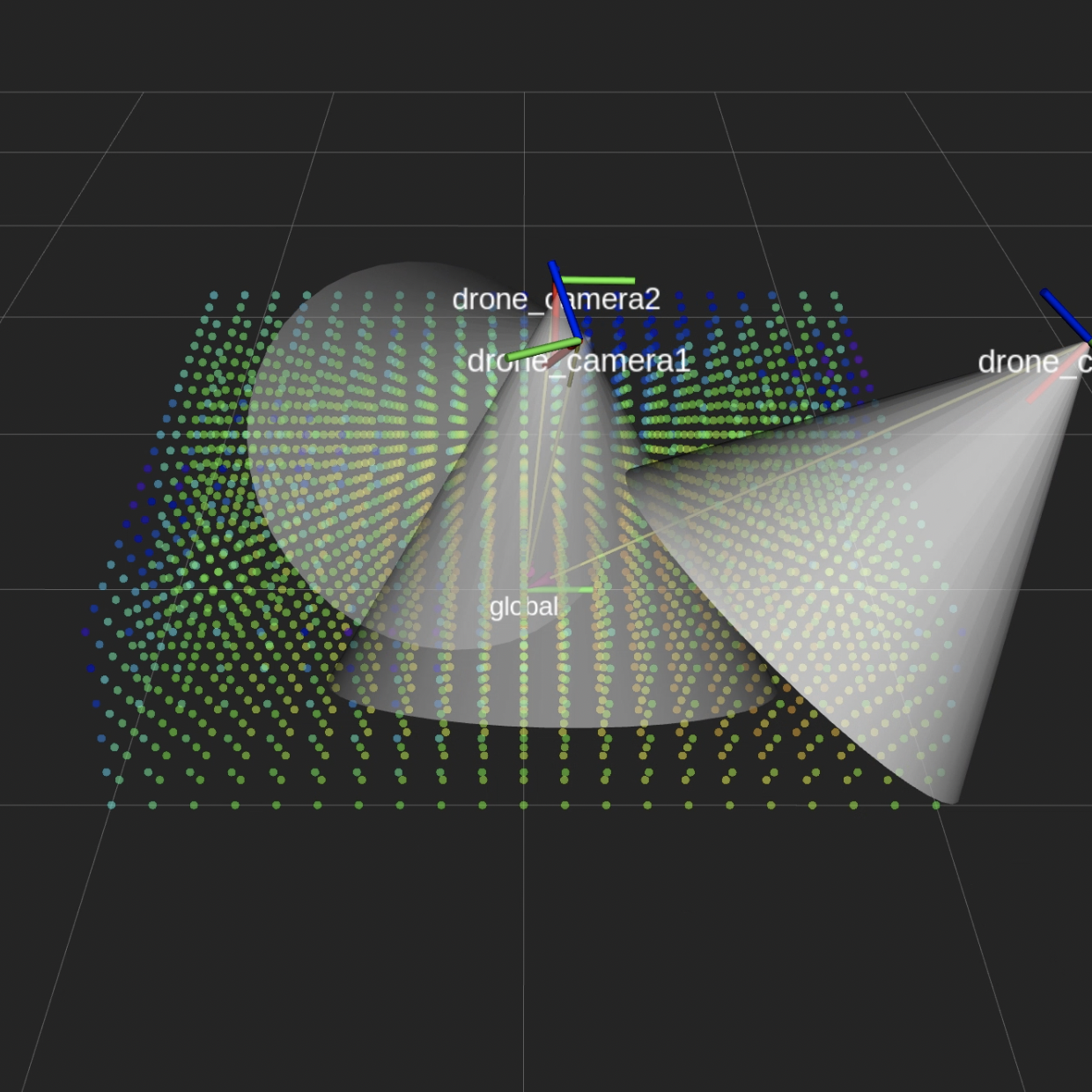}}}\hspace{5pt}
    \subfloat[$t=150s$]{\resizebox*{4cm}{!}{\includegraphics{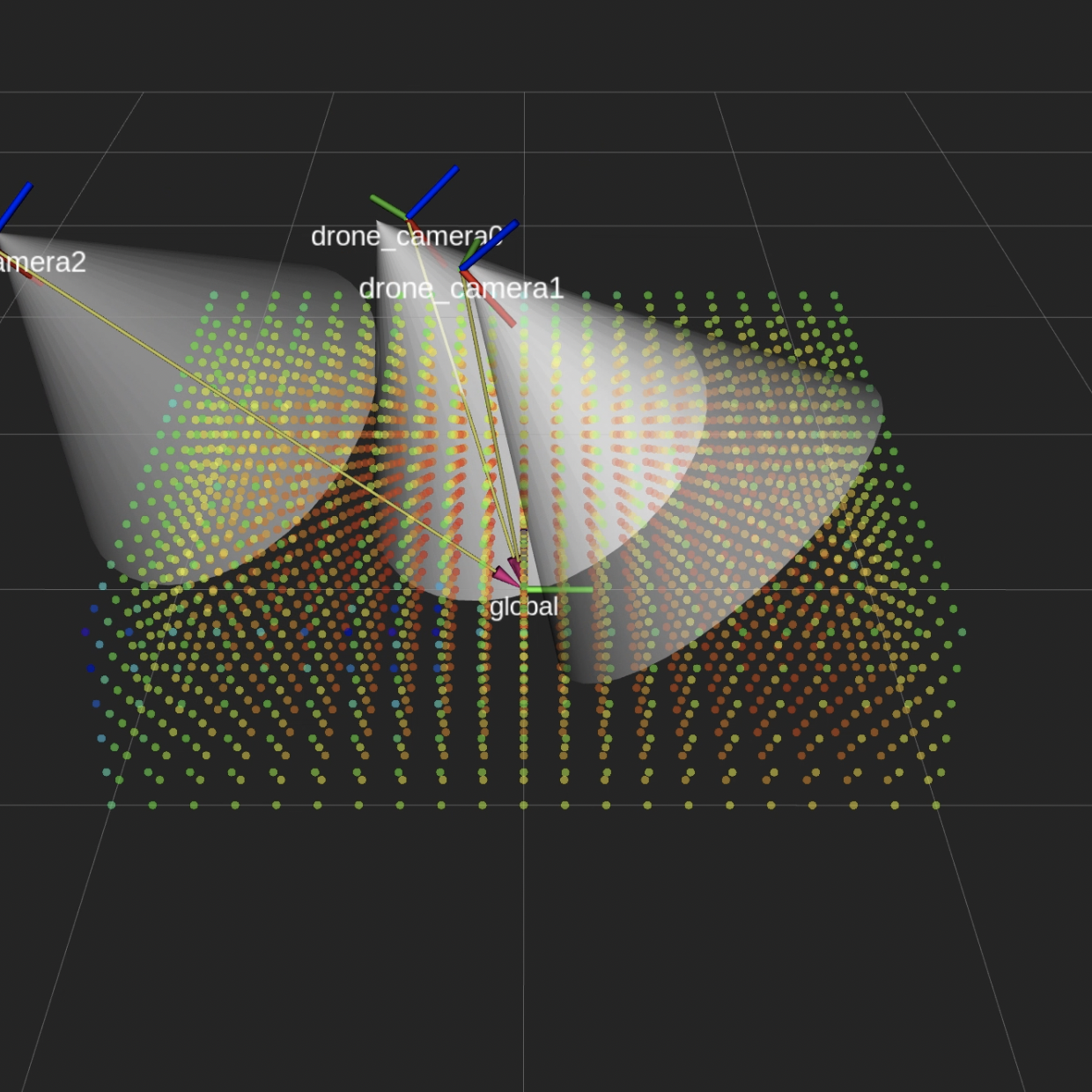}}}\hspace{5pt}
    \subfloat[$t=300s$]{\resizebox*{4cm}{!}{\includegraphics{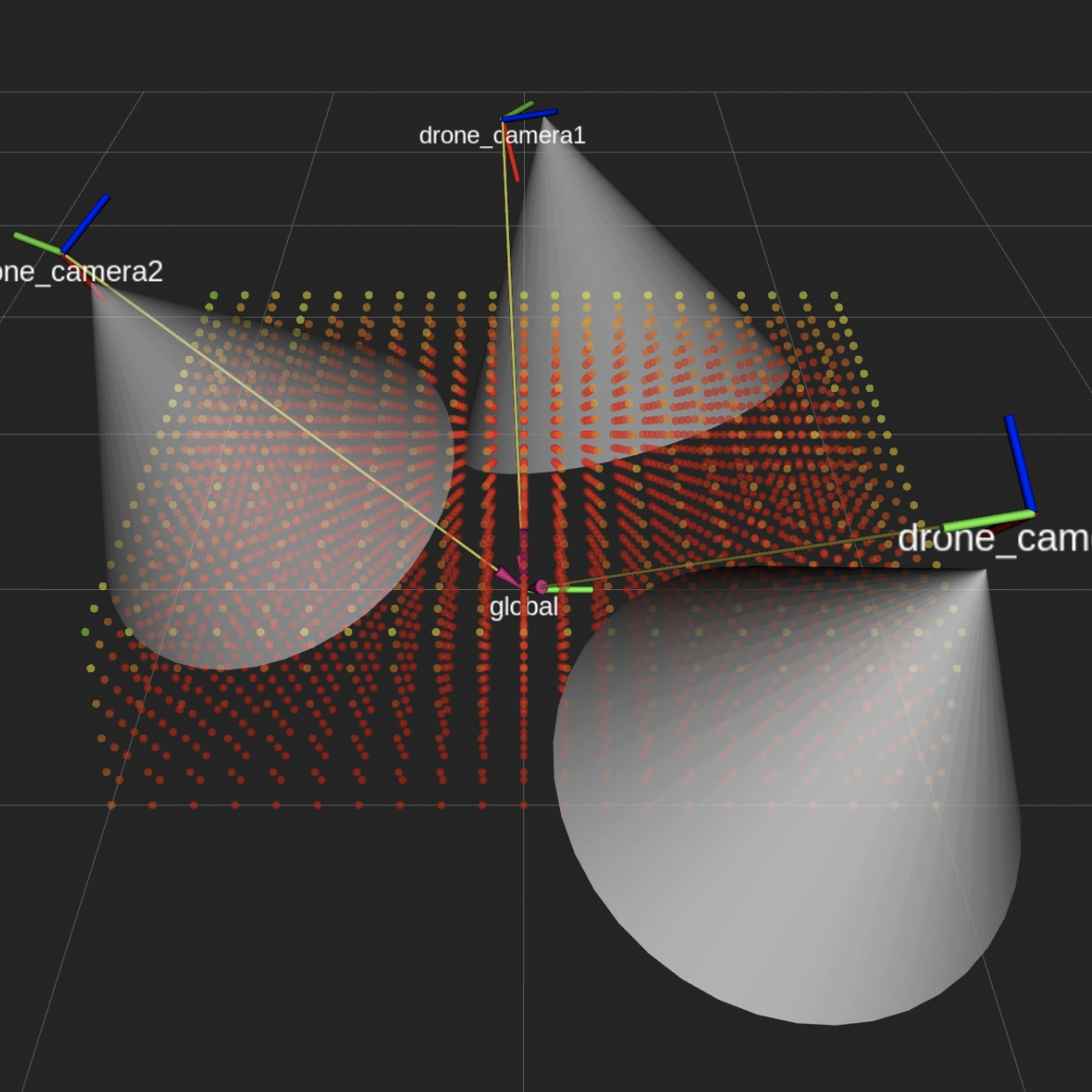}}}\hspace{5pt}
    \caption{Snapshots of the simulation.}
    \label{fig:sim}
    \medskip

    \includegraphics[width = 5.5cm]{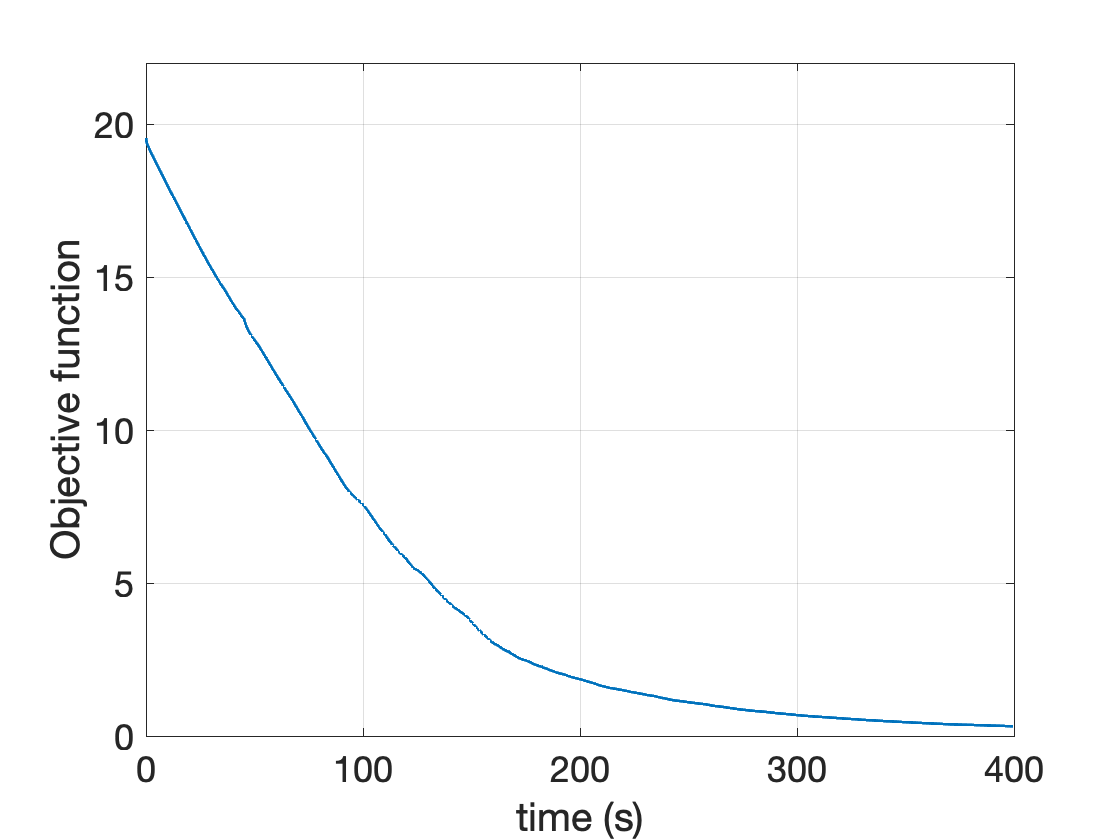}
    \caption{Time evolution of the objective function.}
    \label{fig:obj_func}
\end{figure}

All the drones start the coverage from the same moment $t=0$s, and their positions and camera orientations are presented by a series of screenshots as Fig. \ref{fig:sim}. The value of the importance function is a five-dimensional matrix that is difficult to visualize. Therefore we take the average of the two dimensions $\theta_h$ and $\theta_v$, transform it into a matrix containing $xyz$, and visualize it as a point cloud. The color of the points represents the importance of the position, purple means not yet covered, and red means well covered. The region changed from purple, gradually to red, representing a good completion of the coverage. We observe that drones tend to cover within the task region and keep downwards at the beginning, and start moving out of the task region and turning the gimbal to cover more points at a later stage.

The evolution of the objective function $J$ is shown in Fig. \ref{fig:obj_func}, indicating a linear decreasing trend for the majority of the time, which implies that the decrease of the objective function meets the requirement of the target $\gamma$. After the coverage is completed to a certain extent, the coverage performance of $\gamma$ becomes difficult to achieve and the constraint is violated, so the decline of the objective function gradually slows down.

\begin{figure}[t]
    \centering
    \includegraphics[width = 8cm]{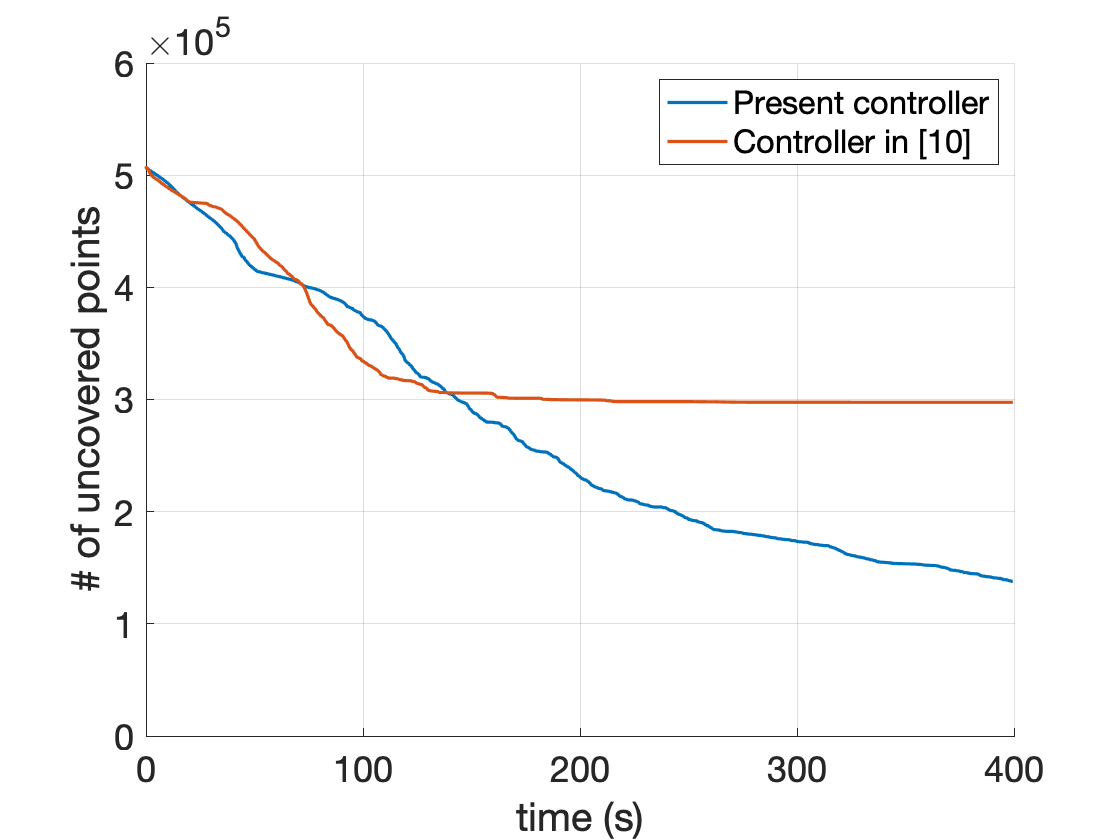}
    \caption{Comparison of angle-aware coverage and with camera rotational motion control}
    \label{fig:comparison}
\end{figure}

Let us next compare the performance of the present controller with that in \cite{shimizu2021angle} without considering the camera orientation control.
Since these controllers employ different objective functions, they cannot be used as metrics of the comparison.
To fairly compare these controllers, we need to prepare another metric that quantifies their performances. In this paper, we revert to the original use case to establish quantifiable metrics related to the actual 3D model quality. In practical applications, drones perceive the environment by taking photos or videos at a limited frequency, rather than through continuous performance functions. We mimic this behavior by attempting to update the drone's coverage of the target region at a finite frequency.
Specifically, we set a shooting rate of images for the drone, representing the number of photos it takes per second. During each shoot, we record which points are covered by the drone. When a point has met the standard of covering, it is marked as ``covered''. This process is similar to updating the objective function, but the difference lies in using Boolean values instead of continuous values. 
By tracking the reduction in the number of uncovered points, we can evaluate how well a controller contributes to the reconstruction of the 3D model. However, we still need to establish criteria to determine whether a point is considered covered. Similar to the performance function, we also consider two aspects. A point is considered covered at a specific moment if and only if both of the following conditions are met: 
\begin{enumerate}
\item The point lies within the field of view of the drone's camera at that moment.
\item The angle between the line connecting the drone to that point and the observed angle of that point differs by less than a certain threshold.
\end{enumerate}

We set the threshold value to $\pi/16$ and recorded the covered points for both of the present controller and \cite{shimizu2021angle} at a burst rate of $5$ Hz. The recorded results of decreased objective functions are presented in Fig. \ref{fig:comparison}.
Since \cite{shimizu2021angle} cannot rotate the camera, the coverage was completed at around $t=130$s and then stopped moving. Due to the limited range of viewing angles, the number of uncovered points remain high even after finishing the coverage. Meanwhile, when we use the present controller with camera orientation control, the number decreases to about 1/2 of \cite{shimizu2021angle}.
This result demonstrates the benefit of controlling the camera orientations as well as the drone positions.

\section{Conclusion}

In this paper, we presented a novel angle-aware coverage with camera rotational motion control. 
To this end, we presented a novel problem formulation including a novel performance function that integrates the camera orientations.
The real-time viability of the present QP-based controller was also demonstrated with the help of JIT and GPU computing.
Moreover, we verified that the present controller achieves a better coverage performance than the original algorithm without camera control \cite{shimizu2021angle}.

Future work should be directed to the hardware experiments and the performance
comparison in terms of the accuracy of the reconstructed 3D structure.

\section*{Acknowledgment}

    This research was partially funded by Japan Society for the Promotion of Science (JSPS) KAKENHI under grant $21K04104$. Additionally, we extend our thanks to Fabrizio Dabbene for his invaluable guidance and to Martina Mammarella for their dedicated contributions, both instrumental in the success of our research.

\section*{Disclosure Statement}

    No potential conflict of interest was reported by the author(s).

\appendix
 
\section{Proof of Theorem \ref{th:qp_form}} \label{app: theorem_1}

The constraint in QP-based controller (\ref{eq:qp_controller}) consists of term $\dot b_{i,I}$, $\dot b_{i,\phi}$, $\alpha_1(b_{i,I})$, and $\alpha_2(b_{i,\phi})$. The term $\dot b_{i,I}$ and $\dot b_{i,\phi}$ can be calculated as follows:

\begin{align}
    \dot b_{i,I} &= \dot I_i = \int_{j\in \mathcal{V}_{i}(p)} \left( \delta \left( \frac{\partial h(p_{i},q_{j})}{\partial p_i} \right) \psi_j + \delta h(p_{i},q_{j}) \dot \psi_j \right) \nonumber \\
     &= \int_{j\in \mathcal{V}_{i}(p)} \left( \delta \left( \frac{\partial h(p_{i},q_{j})}{\partial p_i} \right) \psi_j - \delta^2 h^2(p_{i},q_{j}) \psi_j \right) , \label{eq:b_iI dot}
\end{align}

\begin{align} \label{eq:b_i_phi dot}
    \dot b_{i,\phi} = \left[\begin{array}{c}
        0 \\
        0 \\
        0 \\
        \frac{\partial b_{i,\phi}}{\partial \phi_i}
    \end{array}\right] ^T =\left[\begin{array}{c}
        0 \\
        0 \\
        0 \\
        2 \left ( \phi_{i} - \frac{\phi_{\text{{min}}} + \phi_{\text{{max}}}}{2} \right)
    \end{array}\right] ^ T .
\end{align}

Then, remarking $\alpha_1(b_{i,I}) = a_1 b_{i,I}$ and $\alpha_2(b_{i,\phi}) = a_2 b_{i,\phi}$, we obtain

\begin{align}\label{eq:alpha_biI}
    \alpha_1(b_{i,I}) = a_1 I_i - a_1 \gamma / n = - \frac{a_1 \gamma}{n} +  \left( \int_{j\in \mathcal{V}_{i}(p)} a_1 \delta h(p_{i},q_{j})\psi_{j} \right) ,
\end{align}

\begin{align}\label{eq:alpha_bi_phi}
    \alpha_2(b_{i,\phi}) = a_2 \left ( \left(\phi_{\text{{max}}} - \phi_{\text{{min}}}\right)^2 - \left( \phi_{i} - \frac{\phi_{\text{{min}}} + \phi_{\text{{max}}}}{2} \right)^2 \right) .
\end{align}

By substituting $I_i$ with (\ref{eq:I_i}) and equation (\ref{eq:b_iI dot}),(\ref{eq:alpha_biI}) to constraint (\ref{eq:qp_controller_2}), we can get that $\dot b_{i,I} + \alpha_1(b_{i,I}) = \xi_{1i} u_i + \xi_{2i}$ holds true. Similarly, by substituting equation (\ref{eq:b_i_phi dot}),(\ref{eq:alpha_bi_phi}) to constraint (\ref{eq:qp_controller_3}), we can also get that $\dot b_{i,\phi} + \alpha_2(b_{i,\phi}) = \chi_{1i} u_i + \chi_{2i}$ holds true. Thus, the controller (\ref{eq:qp_controller_final}) is equivalent to (\ref{eq:qp_controller}).

\bibliographystyle{tfnlm}
\bibliography{main}

\end{document}